\documentclass[11pt,a4paper,oneside,reqno]{amsart}
\usepackage[T1]{fontenc}
\usepackage{lmodern,amssymb,graphicx}
\usepackage[a4paper,textwidth=160mm,top=24mm,bottom=25mm,includeheadfoot]{geometry}
\usepackage{array,booktabs,needspace}
\usepackage[hidelinks]{hyperref}
\expandafter\def\expandafter\UrlBreaks\expandafter{\UrlBreaks\do\a\do\b\do\c\do\d\do\e\do\f\do\0\do\1\do\2\do\3\do\4\do\5\do\6\do\7\do\8\do\9}
\hypersetup{pdfauthor={Yasuhito Tsujimura},pdftitle={Triangle-overlap cores and density bounds for Tuza's conjecture},bookmarksnumbered=false}
="[ipaexm.ttf]" at 10.95pt
\providecommand{\tightlist}{\setlength{\itemsep}{2pt}\setlength{\parskip}{0pt}}
\makeatletter
\newcommand{\statementtarget}[1]{\Hy@raisedlink{\hypertarget{#1}{}}}
\renewcommand\subsection{\@startsection{subsection}{2}\z@{1.8ex plus .4ex minus .2ex}{.8ex plus .2ex}{\normalfont\bfseries}}
\newenvironment{papertable}{\par\addvspace{9pt}\noindent\begin{minipage}{\linewidth}\def\@captype{table}}{\end{minipage}\par\addvspace{9pt}}
\newenvironment{paperfigure}{\par\addvspace{9pt}\noindent\begin{minipage}{\linewidth}\def\@captype{figure}}{\end{minipage}\par\addvspace{9pt}}
\def\ps@paper{\def\@oddhead{\hfil\normalfont\scriptsize\scshape Triangle-overlap cores and density bounds\hfil}\let\@evenhead\@oddhead\def\@oddfoot{\hfil\normalfont\small\thepage\hfil}\let\@evenfoot\@oddfoot}
\makeatother
\title[Triangle-overlap cores]{Triangle-overlap cores and density bounds\\for Tuza's conjecture}
\author{Yasuhito Tsujimura}
\date{}
\newcommand{\printauthorinfo}{\normalfont
  Yasuhito Tsujimura\par\nobreak\indent
  Universion Closure Logic\par\nobreak\indent
  Email: \href{mailto:yas@universion-closure-logic.com}{\nolinkurl{yas@universion-closure-logic.com}}\par\nobreak\indent
  ORCID: \href{https://orcid.org/0009-0004-9964-9607}{\nolinkurl{0009-0004-9964-9607}}\par\nobreak\indent
  Website: \url{https://universion-closure-logic.com/}}
\address{\printauthorinfo}
\begin{document}
\begin{abstract}
\small
For a finite simple graph \(G\), let \(\nu(G)\) be the maximum number of
edge-disjoint triangles and let \(\tau(G)\) be the minimum number of
edges meeting every triangle. We prove \(\tau(G)\le2\nu(G)\) when each
edge belongs to at most four triangles, the edges belonging to four
triangles form components with at most two edges, and the maximum
average degree is at most \(25/2\). No upper bound on maximum vertex
degree is imposed. We also localize the density condition. Let \(K\) be
the Hifumi core, defined as the vertex set of the 4-core of the spanning
subgraph formed by edges belonging to at least three triangles. It
suffices to impose the density condition on \(G[K]\), including all
original edges between vertices of \(K\). More generally, when the
four-triangle edges have maximum degree at most two, reductions produce
a subgraph \(R\subseteq G[K]\) with
\(\tau(G)-2\nu(G)\le\tau(R)-2\nu(R)\). Consequently an empty core gives
a density-free sufficient condition. For the first graph class, a
quantitative surplus argument strengthens both the global and localized
density conditions to the strict bound \(25/2+2/56775\). The proof
combines edge reductions, discharging, finite certificates and a weak
König--Egerváry lemma for links of arbitrary order. The general
conjecture and the removal of the density condition from the first graph
class remain unresolved here.
\end{abstract}
\maketitle
\begingroup\small
\textbf{Computational supplement.} The companion verification package
supplies the executable finite checks and an offline reader guide.
Section 6.1 describes its scope and reserved DOI. The package will be
made publicly available after this article is announced on arXiv.
\par\endgroup

\Needspace{5\baselineskip}
\section{1. Results and context}\label{results-and-context}

All graphs are finite, simple and undirected. A \emph{triangle packing}
is a family of pairwise edge-disjoint triangles, and a \emph{triangle
edge cover} is a set of edges meeting every triangle. Their maximum and
minimum sizes are denoted by \(\nu(G)\) and \(\tau(G)\), respectively.
The union of the edges of a maximum packing is a cover, so

\[
\nu(G)\le\tau(G)\le3\nu(G).
\]

Tuza's conjecture asks whether the factor three can always be replaced
by two; see \hyperlink{bib-1}{[1]} for its background and the reduction framework used
below.

This paper addresses two questions: how far a density bound can
guarantee the conjecture under restrictions on triangle overlap, and
where that density bound is actually needed. We first prove Tuza's
inequality for graphs in the class \(\mathcal C\) defined below whose
maximum average degree is at most \(25/2\). Within this class, it is
enough to impose the density bound on the graph induced by a canonical
set of vertices, the Hifumi core. The localization theorem holds in the
larger class \(\mathcal D\); if its core is empty, no density assumption
is needed. These conclusions are stated separately because their overlap
hypotheses differ.

\Needspace{14\baselineskip}
For \(uv\in E(G)\), define its \emph{codegree} \[
b_G(uv)=|N_G(u)\cap N_G(v)|,
\] the number of triangles containing \(uv\). Write
\(b(G)=\max_{e\in E(G)}b_G(e)\), with \(b(G)=0\) when \(G\) is edgeless.
Let \(H_i(G)\) be the spanning subgraph with edge set
\(\{e:b_G(e)=i\}\), and let \(H_{\ge3}(G)\) have edge set
\(\{e:b_G(e)\ge3\}\). The maximum average degree is \[
\operatorname{mad}(G)=
\max_{\substack{J\subseteq G\\|V(J)|>0}}
\frac{2|E(J)|}{|V(J)|};
\qquad \operatorname{mad}(\varnothing)=0.
\tag{1.1}\label{eq:1.1}
\] In particular, a bound on \(\operatorname{mad}(G)\) is a bound on
every subgraph, not merely on the average degree of \(G\).

\Needspace{9\baselineskip}
We use two hereditary classes:

\begin{itemize}
\tightlist
\item
  \(\mathcal C\): \(b(G)\le4\), and each nontrivial component of
  \(H_4(G)\) is \(K_2\) or \(P_3\), where \(P_3\) has three vertices and
  two edges;
\item
  \(\mathcal D\): \(b(G)\le4\) and \(\Delta(H_4(G))\le2\).
\end{itemize}

Thus \(\mathcal C\subseteq\mathcal D\). In \(\mathcal D\), the
nontrivial components of \(H_4\) may be paths of arbitrary length or
cycles. Neither class imposes an upper bound on \(\Delta(G)\).

In the common subcase \(b(G)\le3\), our density bound replaces \(45/4\)
in the supplementary statement \hyperlink{bib-3}{[3]} by \(25/2\). The comparison is
confined to this subcase; \hyperref[relation-to-earlier-results]{Section~1.1} distinguishes the hypotheses and
attribution of the related results.

\par\addvspace{7pt}\Needspace{5\baselineskip}\noindent

\textbf{\statementtarget{statement-theorem-1-1}Theorem 1.1 (density bound).} If \(G\in\mathcal C\) and
\(\operatorname{mad}(G)\le25/2\), then \(\tau(G)\le2\nu(G)\).

The endpoint \(25/2\) is included. To state the second result, put \[
\Phi(G)=\tau(G)-2\nu(G),\qquad
K=V\bigl(\operatorname{core}_4(H_{\ge3}(G))\bigr).
\tag{1.2}\label{eq:1.2}
\] Here the \(k\)-core is the maximal induced subgraph of minimum degree
at least \(k\), obtained by repeatedly removing vertices of current
degree below \(k\). We call the vertex set \(K\) the \emph{Hifumi
triangle-overlap core}, or \emph{Hifumi core} for short. It is computed
from the initial graph. The induced graph \(G[K]\) retains \textbf{all}
original edges between these vertices, including those of codegree less
than three.

\par\addvspace{7pt}\Needspace{5\baselineskip}\noindent

\textbf{\statementtarget{statement-theorem-1-2}Theorem 1.2 (Hifumi core localization).} For every
\(G\in\mathcal D\), a sequence of edge reductions and deletions of
isolated vertices yields a subgraph \(R\subseteq G[K]\) such that \[
\Phi(G)\le\Phi(R).
\tag{1.3}\label{eq:1.3}
\] The reductions explicitly extend packings and covers from the
surviving graph to the preceding graphs.

The survivor may omit edges inside \(K\); it need not be \(G[K]\)
itself. This distinction is essential when applying a sufficient
condition to the core.

\begin{figure}[!htbp]
\centering
\includegraphics[width=\linewidth]{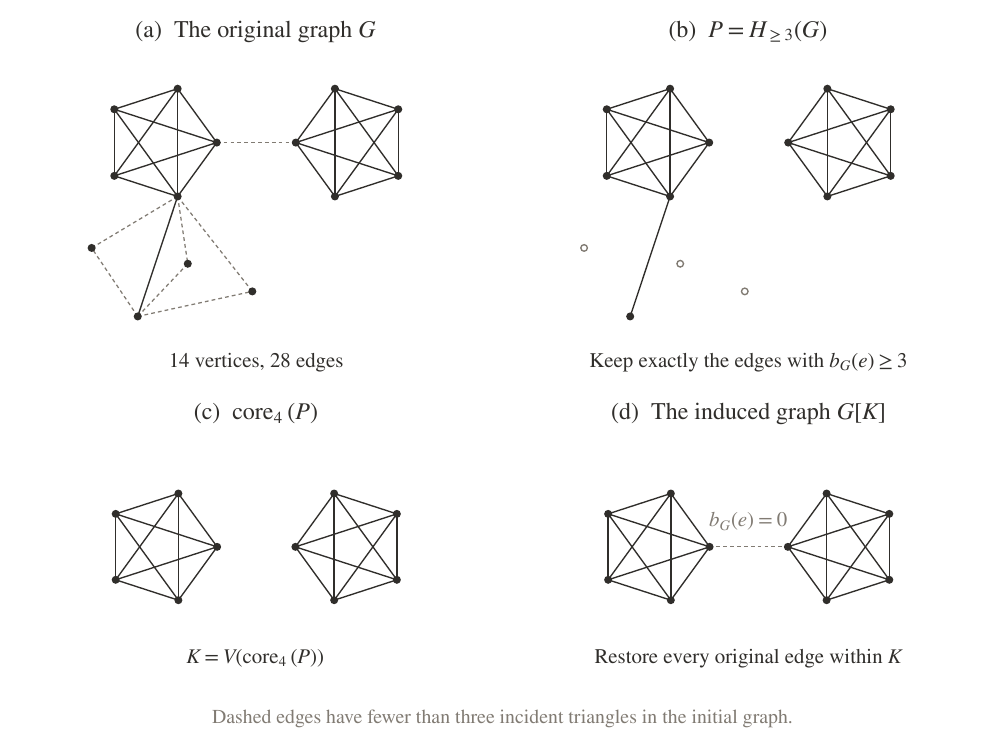}
\caption{The Hifumi core and the induced graph on its vertices. The graph in (a)
consists of two copies of \(K_5\), joined by an edge, with a three-page
book attached at one clique vertex. Panel (b) retains exactly the edges
of \(P=H_{\ge3}(G)\); isolated vertices are shown as open circles.
Removing vertices of current degree below four gives (c), whose ten
vertices form \(K\). Panel (d) restores all original edges with
endpoints in \(K\). In particular, the joining edge, of initial codegree
zero, belongs to \(G[K]\). Dashed edges in (a) and (d) have initial
codegree less than three. Only marked points are vertices; crossings are
not vertices. The figure illustrates the definition of \(K\), not the
survivor \(R\) of Theorem 1.2.}
\label{fig-hifumi-core}
\hypertarget{fig-hifumi-core}{}
\end{figure}

\hyperref[fig-hifumi-core]{Figure 1} illustrates why the vertex set
\(K\) and the induced graph \(G[K]\) must be distinguished from the
high-codegree subgraph used to compute the core.

\par\addvspace{7pt}\Needspace{5\baselineskip}\noindent

\textbf{\statementtarget{statement-corollary-1-3}Corollary 1.3.} Each of the following implies
\(\tau(G)\le2\nu(G)\):

\begin{enumerate}
\def\labelenumi{\arabic{enumi}.}
\tightlist
\item
  \(G\in\mathcal C\) and \(\operatorname{mad}(G[K])\le25/2\);
\item
  \(G\in\mathcal D\) and \(H_{\ge3}(G)\) is 3-degenerate, equivalently
  \(K=\varnothing\).
\end{enumerate}

Indeed, in the first case every subgraph of \(G[K]\) satisfies the
hypotheses of Theorem 1.1. In the second case \(R\) is empty. Equation
\hyperref[eq:1.3]{(1.3)} gives the conclusion. The first condition imposes no density bound
outside \(K\); the second imposes no density bound at all.

\Needspace{16\baselineskip}
Set \(\eta=2/56775\).

\par\addvspace{7pt}\Needspace{5\baselineskip}\noindent

\textbf{\statementtarget{statement-theorem-1-4}Theorem 1.4 (quantitative gap).} A nonempty graph in
\(\mathcal C\), with no isolated vertices and admitting no edge
reduction of \hyperlink{statement-definition-2-1}{Definition~2.1}, has average degree at least \(25/2+\eta\).
Consequently \hyperlink{statement-theorem-1-1}{Theorem~1.1} and \hyperlink{statement-corollary-1-3}{Corollary~1.3}(1) hold with their density
hypothesis replaced by \[
\operatorname{mad}(G)<\frac{25}{2}+\eta
\quad\text{or, respectively,}\quad
\operatorname{mad}(G[K])<\frac{25}{2}+\eta.
\tag{1.4}\label{eq:1.4}
\] Since \(1/30000<\eta\), the corresponding closed bound
\(25/2+1/30000\) is also sufficient. We claim no optimality for these
thresholds. The gap expresses that the excluded equality structure
forces a uniform positive average surplus.

\hyperref[tab-results]{Table 1} summarizes the main results. The full
hypotheses and the quantitative refinements are stated above.

\Needspace{10\baselineskip}
\begin{papertable}
\caption{Main results and their proof locations. Refined density thresholds are
given in Theorem 1.4.}
\label{tab-results}
\hypertarget{tab-results}{}
\vspace{4pt}\centering\small\renewcommand{\arraystretch}{1.22}
\begin{tabular}{@{}>{\raggedright\arraybackslash}p{25.670mm}@{\hspace{3mm}}>{\raggedright\arraybackslash}p{39.260mm}@{\hspace{3mm}}>{\raggedright\arraybackslash}p{57.380mm}@{\hspace{3mm}}>{\raggedright\arraybackslash}p{28.690mm}@{}}
\toprule
\textbf{Result} & \textbf{Hypotheses} & \textbf{Conclusion} & \textbf{Proof} \\
\midrule
\hyperlink{statement-theorem-1-1}{Theorem~1.1} & \(G\in\mathcal C\), \(\operatorname{mad}(G)\le25/2\) & \(\tau(G)\le2\nu(G)\) & \hyperref[discharging-and-the-density-bound]{Section~3}; Appendices A--B \\
\hyperlink{statement-theorem-1-2}{Theorem~1.2} & \(G\in\mathcal D\) & Reductions give \(R\subseteq G[K]\) with \(\Phi(G)\le\Phi(R)\) & \hyperref[localization-to-the-hifumi-core]{Section~4}; \hyperref[appendix-c.-wke-witnesses-for-links-of-arbitrary-order]{Appendix~C} \\
\hyperlink{statement-corollary-1-3}{Corollary~1.3}(1) & \(G\in\mathcal C\), \(\operatorname{mad}(G[K])\le25/2\) & \(\tau(G)\le2\nu(G)\); no density condition outside \(K\) & Theorems 1.1--1.2 \\
\hyperlink{statement-corollary-1-3}{Corollary~1.3}(2) & \(G\in\mathcal D\), \(K=\varnothing\) & \(\tau(G)\le2\nu(G)\); no density condition & \hyperlink{statement-theorem-1-2}{Theorem~1.2} \\
\hyperlink{statement-theorem-1-4}{Theorem~1.4} & Nonempty reduction-free \(G\in\mathcal C\), with no isolated vertices & Average degree at least \(25/2+\eta\) & \hyperref[a-uniform-positive-surplus]{Section~5}; \hyperref[appendix-d.-the-quantitative-surplus-count]{Appendix~D} \\
\bottomrule
\end{tabular}
\end{papertable}

\Needspace{5\baselineskip}
\subsection{1.1. Relation to earlier
results}\label{relation-to-earlier-results}

Puleo \hyperlink{bib-1}{[1, Theorem 1.2]} established the conjecture for maximum
average degree strictly less than seven, without our overlap
restrictions. We use his edge-reduction and weak König--Egerváry
framework \hyperlink{bib-1}{[1, Definitions 4.1--4.2 and Lemmas 4.3--4.4]}. Our larger
density threshold has additional structural hypotheses and is not an
unconditional improvement of that theorem.

Kandasamy \hyperlink{bib-2}{[2, Theorem 1.1]} proves the conjecture for all graphs of
maximum degree at most ten. A supplementary document in the verification
archive linked by that paper \hyperlink{bib-3}{[3]} states the sufficient conditions
\(b(G)\le3\) and \(\operatorname{mad}(G)\le45/4\). In the common subcase
\(b(G)\le3\), \hyperlink{statement-theorem-1-1}{Theorem~1.1} replaces \(45/4\) by \(25/2\). The
four-triangle-edge case and localization concern additional structure.
Removing a maximum-degree cap does not mean that our hypotheses include
every graph covered by \hyperlink{bib-2}{[2]}.

Liu and Zeng \hyperlink{bib-4}{[4, Theorem 1.3]} prove the conjecture for every fixed
real \(c<2\), for sufficiently large \(n\)-vertex graphs with at least
\(n^2/4-cn\) edges. This is a different density regime. Munaro \hyperlink{bib-5}{[5,
Theorem 3]} proves the stronger ratio \(\tau\le3\nu/2\) when each
triangle shares at most two of its edges with other triangles and each
edge belongs to at most four triangles. The shared-edge condition is
additional; that result is not a theorem for every graph with
\(b(G)\le4\).

\Needspace{5\baselineskip}
\subsection{1.2. Structure of the proof}\label{structure-of-the-proof}

An edge reduction preserves a counterexample: if a graph violates Tuza's
inequality, so does the smaller graph left by the reduction (\hyperlink{statement-lemma-2-2}{Lemma~2.2}).
The two main arguments then separate.

For the density theorem, apply \hyperlink{statement-lemma-2-5}{Lemma~2.5} and the discharging rules of
\hyperref[discharging-and-the-density-bound]{Section~3} to a nonempty reduction-free graph in \(\mathcal C\) without
isolated vertices. Every vertex retains charge at least \(25/2\), with
strict surplus at degree at least fourteen. \hyperref[appendix-b.-the-degree-twelvethirteen-endpoint]{Appendix~B} excludes the
degree-twelve/thirteen equality case, completing \hyperlink{statement-theorem-1-1}{Theorem~1.1} in Section
3.4.

For localization, \hyperlink{statement-lemma-4-1}{Lemma~4.1} and \hyperref[appendix-c.-wke-witnesses-for-links-of-arbitrary-order]{Appendix~C} provide WKE witnesses for
links of arbitrary order. \hyperref[localization-to-the-hifumi-core]{Section~4} uses them to reduce outside the
initial Hifumi core, proving \hyperlink{statement-theorem-1-2}{Theorem~1.2} independently of the density
theorem. \hyperlink{statement-corollary-1-3}{Corollary~1.3} combines localization with hereditary sufficient
conditions.

\hyperref[a-uniform-positive-surplus]{Section~5} and \hyperref[appendix-d.-the-quantitative-surplus-count]{Appendix~D} give Theorem 1.4. \hyperref[verification-scope-and-a-remaining-question]{Section~6} specifies the
verification boundary and remaining question. Appendices A--E form part
of the proof; \hyperref[appendix-e.-an-ordinary-core-sufficient-condition]{Appendix~E} adds a sufficient condition using ordinary
graph cores.

\Needspace{5\baselineskip}
\section{2. Reductions and local
structure}\label{reductions-and-local-structure}

\Needspace{5\baselineskip}
\subsection{2.1. Edge reductions}\label{edge-reductions}

\par\addvspace{7pt}\Needspace{5\baselineskip}\noindent

\textbf{\statementtarget{statement-definition-2-1}Definition 2.1.} An \emph{edge reduction} consists of a nonempty
set \(F\subseteq E(G)\), an edge-disjoint triangle family
\(\mathcal S\), and \(X\subseteq E(G)\) satisfying \[
\begin{split}
&X\text{ meets every triangle containing an edge of }F,\\
&E(\mathcal S)\setminus F\subseteq X,\qquad |X|\le2|\mathcal S|.
\end{split}
\tag{2.1}\label{eq:2.1}
\] This is the edge-reducibility condition of \hyperlink{bib-1}{[1, Definition 4.2]},
with the deleted edge set denoted here by \(F\). No graph-minor
contraction is involved.

\par\addvspace{7pt}\Needspace{5\baselineskip}\noindent

\textbf{\statementtarget{statement-lemma-2-2}Lemma 2.2 (reduction inequality).} For \(G'=G-(F\cup X)\), \[
\nu(G)\ge\nu(G')+|\mathcal S|,
\qquad \tau(G)\le\tau(G')+|X|,
\qquad \Phi(G)\le\Phi(G').
\tag{2.2}\label{eq:2.2}
\]

\emph{Proof.} Every edge of a triangle in \(\mathcal S\) has been
removed from \(G'\). Thus any packing of \(G'\) can be augmented by
\(\mathcal S\). A triangle of \(G\) absent from \(G'\) uses an edge of
\(F\cup X\); it is met by \(X\), directly or by \hyperref[eq:2.1]{(2.1)}. A cover of
\(G'\), together with \(X\), covers \(G\). Subtraction gives the last
inequality. This spells out the extension in \hyperlink{bib-1}{[1, Lemma 4.3]}. \hfill\(\square\)

\begin{figure}[!htbp]
\centering
\includegraphics[width=\linewidth]{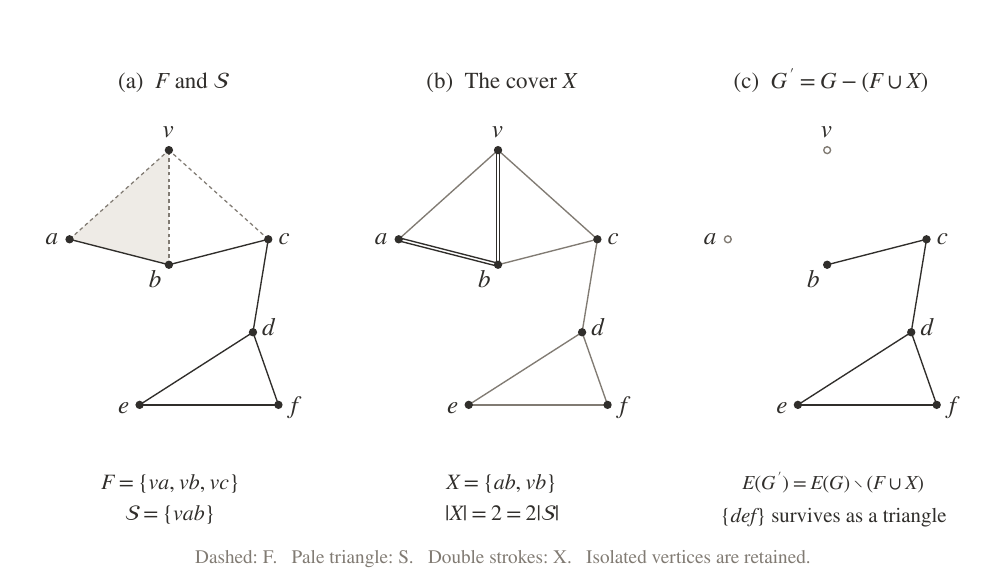}
\caption{A literal edge reduction. In (a), dashed edges form \(F=\{va,vb,vc\}\),
and the shaded triangle is the sole member of \(\mathcal S=\{vab\}\). In
(b), double strokes mark the single edges in \(X=\{ab,vb\}\); they do
not represent parallel edges. Every triangle using an edge of \(F\)
meets \(X\), \(E(\mathcal S)\setminus F=\{ab\}\subseteq X\), and
\(|X|=2|\mathcal S|\). Panel (c) removes \(F\cup X\), retaining the
isolated vertices \(v,a\) as open circles. The surviving triangle
\(def\) can be packed together with \(vab\) in the original graph. No
vertices are identified.}
\label{fig-edge-reduction}
\hypertarget{fig-edge-reduction}{}
\end{figure}

\hyperref[fig-edge-reduction]{Figure 2} shows a reduction with
\(F\cap X\ne\varnothing\) and \(X\setminus F\ne\varnothing\), making the
deletion of \(F\cup X\) explicit.

Isolated vertices do not affect \(\tau,\nu,\Phi\). In structural lemmas,
a \emph{reduction-free} graph admits no reduction of \hyperlink{statement-definition-2-1}{Definition~2.1}, and
isolated vertices are removed. An input to our main theorems need not
already be reduction-free.

\par\addvspace{7pt}\Needspace{5\baselineskip}\noindent

\textbf{\statementtarget{statement-lemma-2-3}Lemma 2.3 (WKE reduction).} The \emph{link} of \(v\) is
\(L_v=G[N_G(v)]\). Suppose a nonempty union \(L'\) of components of
\(L_v\) has a matching \(M\) and a vertex cover \(Q\) of \(L'-M\), with
\(|Q|\le|M|\). Then \(G\) has an edge reduction deleting all spokes from
\(v\) to \(V(L')\).

\emph{Proof.} Such a pair \((M,Q)\) is a weak König--Egerváry, or WKE,
witness \hyperlink{bib-1}{[1, Definition 4.1]}. Set \[
F=\{vu:u\in V(L')\},\quad
\mathcal S=\{vxy:xy\in M\},\quad
X=M\cup\{vq:q\in Q\}.
\tag{2.3}\label{eq:2.3}
\] Matching edges give edge-disjoint triangles. Every triangle using
\(F\) is determined by an edge of \(L'\), either in \(M\) or met by
\(Q\). The remaining conditions follow immediately. This is [1, Lemma
4.4]. \hfill\(\square\)

\par\addvspace{7pt}\Needspace{5\baselineskip}\noindent

\textbf{\statementtarget{statement-lemma-2-4}Lemma 2.4 (heredity and elementary deletions).} The classes
\(\mathcal C,\mathcal D\) are closed under taking subgraphs. In a
reduction-free graph of either class every edge has codegree at least
two, and every link has minimum degree at least two, maximum degree at
most four and at most two vertices of degree four.

\emph{Proof.} Codegrees do not increase on deletion. Under \(b(G)\le4\),
an edge of codegree four in a subgraph already had codegree four in
\(G\); hence \(H_4\) only loses edges. Both restrictions on \(H_4\) are
preserved. A codegree-zero edge is removed with
\(\mathcal S=X=\varnothing\). If \(ab\) lies only in triangle \(abc\),
take \(F=\{ab\}\), \(\mathcal S=\{abc\}\), \(X=\{ac,bc\}\). Finally
\(d_{L_v}(u)=b_G(vu)\) and \(\Delta(H_4)\le2\). \hfill\(\square\)

\Needspace{5\baselineskip}
\subsection{2.2. Structural premises for
discharging}\label{structural-premises-for-discharging}

Call a vertex \emph{high} if its degree in the current graph is at least
fourteen. A label in a link identifies an actual neighbor of its center.
For degree-nine and degree-eleven centers we specify \emph{designated
neighbors} from which charge will be received. The designation is local
to that center.

\par\addvspace{7pt}\Needspace{5\baselineskip}\noindent

\noindent\begin{minipage}{\linewidth}
\textbf{\statementtarget{statement-lemma-2-5}Lemma 2.5 (local structural facts).} In a nonempty
reduction-free \(G\in\mathcal C\):

\begin{enumerate}
\def\labelenumi{(\alph{enumi})}
\item
  Every vertex has degree at least nine. A link of order at most
  seventeen is connected.
\item
  A degree-nine link is one of \(N_1,N_2,N_3,N_4\), defined in Appendix
  A. The first three have one incident codegree-four edge and eight
  incident codegree-three edges; all nine neighbors are designated. In
  the fourth, eight neighbors are designated: two are joined to the
  center by codegree-four edges and six by codegree-three edges. The
  remaining neighbor is joined by a codegree-two edge and is
  undesignated. Every designated neighbor has degree at least fifteen.
  All four common neighbors along a designated codegree-four edge are
  high. Along a codegree-three edge in the first three types, all three
  are high; in \(N_4\) at least two are high, and only labels 2 and 3
  may have exactly two.
\item
  A degree-ten link is \(T_P\) or \(T_R\) in Appendix A. All neighbors
  have degree at least fourteen, and all incident edges have codegree
  three. A degree-fourteen neighbor is possible only for \(T_R\), at
  labels 4 or 9; its link is then \(J_{14}\), with the degree-ten center
  at label 1 or 6.
\item
  Each degree-eleven vertex can designate nine neighbors of degree at
  least fourteen. At each eligible designation the link degree is three
  or four and at least two common neighbors are high. The choice can be
  made independently at each degree-eleven center.
\item
  Each degree-twelve vertex has at most two neighbors of degree at most
  twelve. Each degree-thirteen vertex has at most ten degree-twelve
  neighbors.
\item
  No two vertices of degrees nine or ten are adjacent. In particular the
  designated degree-nine neighbors and degree-ten neighbors of any
  vertex form an independent set in its link.
\end{enumerate}
\end{minipage}\par

\Needspace{6\baselineskip}
\hyperref[appendix-a.-finite-domains-and-local-graph-data]{Appendix~A} specifies the finite domains, certificate conditions and
shared-center facts used at equality. These are finite verification
premises, not extra hypotheses on the input graph. Their application
uses only the implication from an actual graph to the finite domain. A
retained finite configuration need not extend to a graph.

For (a), the broader small-component premise is that every connected
graph of order at most eight, maximum degree at most four and at most
two vertices of degree four has a WKE or direct cone reduction. The cone
over a link component consists of that component together with its
center and all spokes joining the center to the component. A direct cone
certificate satisfies \hyperref[eq:2.1]{(2.1)} with \(F\) the spokes. Every component of a
reduction-free link must have order at least nine; two components would
require eighteen vertices. This premise also applies in \(\mathcal D\),
as used in Appendix E.

\Needspace{5\baselineskip}
\section{3. Discharging and the density
bound}\label{discharging-and-the-density-bound}

Work in a nonempty reduction-free \(G\in\mathcal C\). Initially assign
charge \(d_G(v)\) to each vertex. Choose the designated neighbors of
\hyperlink{statement-lemma-2-5}{Lemma~2.5} once, and use the transfers in \hyperref[tab-transfers]{Table
2}.

\Needspace{10\baselineskip}
\begin{papertable}
\caption{Charge received by low-degree vertices. Designations are fixed as in
Lemma 2.5.}
\label{tab-transfers}
\hypertarget{tab-transfers}{}
\vspace{4pt}\centering\small\renewcommand{\arraystretch}{1.22}
\begin{tabular}{@{}>{\raggedright\arraybackslash}p{32.340mm}@{\hspace{3mm}}>{\raggedright\arraybackslash}p{98.560mm}@{\hspace{3mm}}>{\raggedright\arraybackslash}p{23.100mm}@{}}
\toprule
\textbf{Recipient} & \textbf{Source} & \textbf{Charge} \\
\midrule
degree nine & designated neighbor along a codegree-four edge & \(7/10\) \\
degree nine & designated neighbor along a codegree-three edge & \(7/20\) \\
degree ten & each neighbor & \(1/4\) \\
degree eleven & each of its nine designated neighbors & \(1/6\) \\
degree twelve & each neighbor of degree at least thirteen & \(1/20\) \\
\bottomrule
\end{tabular}
\end{papertable}

Let \(\operatorname{ch}(v)\) be the final charge. The total is
conserved. Vertices of degrees nine, ten and eleven finish with
\(25/2\). A vertex of degree twelve has at least ten neighbors of degree
at least thirteen and finishes with at least \(25/2\). A vertex of
degree thirteen pays at most \(10/20\) and also finishes with at least
\(25/2\).

\Needspace{5\baselineskip}
\Needspace{12\baselineskip}
\subsection{3.1. Bounds at a high vertex}\label{bounds-at-a-high-vertex}

Fix \(h\) of degree \(d\ge14\). Let \(a,b,c,e\) count, respectively, its
degree-nine neighbors designating \(h\), degree-ten neighbors,
degree-eleven neighbors designating \(h\), and degree-twelve neighbors.

Among the high neighbors, let \(k\) count the vertices and \(z\) the
codegree-four edges joining them to \(h\).

Within the first group, let \(q\) count the codegree-four edges to
\(h\). Let \(t\) count the designated degree-nine neighbors joined to
\(h\) by a codegree-three edge with only two high common neighbors.
Other low neighbors may be undesignated and cost nothing.

Counting incidences with high vertices in the link gives \[
\begin{gathered}
q+z\le2,\quad0\le t\le a-q,\quad k\le d-a-b-c-e,\\
3a+q-t+3b+2c\le3k+z.
\end{gathered}
\tag{3.1}\label{eq:3.1}
\] The left side counts the required high common neighbors from Lemma
2.5. A vertex of degree ten contributes three, since all its neighbors
are high. Each high link vertex has degree at most three with at most
\(z\) extra units in total.

The \(a+b\) vertices of degree nine or ten are independent. Their
link-degree sum is \(3(a+b)+q\), while the complementary degree sum is
at most \(3(d-a-b)+(2-q)\). Hence \[
6(a+b)+2q\le3d+2.
\tag{3.2}\label{eq:3.2}
\] The outgoing charge is \[
\operatorname{out}(h)=\frac{7(a+q)}{20}+\frac b4+\frac c6+\frac e{20}.
\tag{3.3}\label{eq:3.3}
\] Using \(t\le a-q\) and the upper bound on \(k\) in \hyperref[eq:3.1]{(3.1)} gives \[
5a+6b+5c+3e\le3d+z-2q.
\] Multiplication by \(7/100\) bounds all coefficients in \hyperref[eq:3.3]{(3.3)} apart
from the separate \(7q/20\). Thus \[
\operatorname{out}(h)\le\frac{21d+7z+21q}{100}
\le\frac{21d+42}{100}.
\tag{3.4}\label{eq:3.4}
\] For every \(d\ge17\), this leaves charge strictly above \(25/2\). At
\(d=16\), enumeration of the nonnegative integer tuples satisfying
\hyperref[eq:3.1]{(3.1)}--\hyperref[eq:3.2]{(3.2)} gives \(\operatorname{out}(h)\le101/30<7/2\). This is a
bounded integer calculation with an explicit domain, not an enumeration
of order-sixteen graphs.

\Needspace{5\baselineskip}
\subsection{3.2. Degree fourteen}\label{degree-fourteen}

Here \(a=q=0\), and the rooted-pair premises give \(b\le2\). Together
with \[
3b+2c\le3k+z,\quad b+c+e+k\le14,\quad z\le2,
\] the integer bounds give \(\operatorname{out}(h)\le3/2\), with
equality only at \[
(b,c,e,k)=(2,6,0,6).
\tag{3.5}\label{eq:3.5}
\] Both degree-ten neighbors would occupy labels 1 and 6 in \(J_{14}\).
The labels \(\{0,7,12,13\}\) adjacent to them must all be high. The
eight remaining labels form \(\{2,3,8,10\}\) and \(\{4,5,9,11\}\), each
inducing \(K_4\) minus an edge. Each set contains at most two designated
degree-eleven neighbors: if it contained three, at least one of its two
internal degree-three vertices would be designated but have at most one
high neighbor, contrary to \hyperlink{statement-lemma-2-5}{Lemma~2.5}(d). Thus \(c\le4\), contradicting
\hyperref[eq:3.5]{(3.5)}. Excluding this profile from the integer domain gives \[
\operatorname{out}(h)\le\frac{22}{15},\qquad
\operatorname{ch}(h)\ge\frac{188}{15}=\frac{25}{2}+\frac1{30}.
\tag{3.6}\label{eq:3.6}
\]

\Needspace{5\baselineskip}
\subsection{3.3. Degree fifteen}\label{degree-fifteen}

The finite obligation bounds outgoing charge for every link and
compatible assignment of neighboring degrees. Generate all connected
order-fifteen links of minimum degree two, maximum degree four and at
most two degree-four vertices. A cone already violating \(\mathcal C\)
is excluded. Other exclusions require a literal WKE or edge-reduction
certificate. No case is excluded because an optimization search failed.

For each remaining link \(L\), enumerate all \(H\subseteq V(L)\),
representing its high ambient neighbors. Interfaces with the complete
degree-nine, ten and eleven templates preserve adjacency among common
neighbors and all forced high neighbors, under every common-neighbor
ordering. These are necessary tests; allowing extra assignments only
weakens the upper bound.

For an actual graph, define \[
I_*=\{v\in N(h):d_G(v)=10\}
\;\cup\;\{v\in N(h):d_G(v)=9,\ v\text{ designates }h\}.
\] This set is independent by \hyperlink{statement-lemma-2-5}{Lemma~2.5}(f). An undesignated degree-nine
neighbor is not in \(I_*\); in particular, the codegree-two position of
\(N_4\) receives no payment as a degree-nine vertex.

Multiplying the five transfer amounts in Table 2 by 240 gives the
integer costs 168, 84, 60, 40 and 12, respectively. For \(v\notin H\),
let \(B_v=40\) if a compatible contact with a degree-eleven vertex
designating the center is possible, and \(B_v=12\) otherwise. This
bounds the payment to every nonhigh neighbor outside \(I_*\), including
degree-twelve vertices and undesignated neighbors.

\Needspace{7\baselineskip}
Exclusions are justified for each rooted contact and identification. The
root exclusions used here remove contacts with degree-nine vertices only
after all eligible degree-nine templates and identifications at that
root have been excluded; they do not remove a compatible contact with a
degree-ten or degree-eleven vertex.

Let \(A_v\) be the maximum of the costs \(168,84,60\) over surviving
eligible contacts with degree-nine or degree-ten vertices at \(v\). Omit
\(v\) from the candidate set only if no such contact survives. For each
independent set \(I\subseteq V(L)\setminus H\) of candidates, the bound
is \[
\sum_{v\notin H}B_v+
\sum_{v\in I}\max\{0,A_v-B_v\}.
\tag{3.7}\label{eq:3.7}
\] The actual set \(I_*\) is one such independent set. Payments to its
members are bounded by \(A_v\); all other nonhigh payments are bounded
by \(B_v\). Hence the total outgoing charge is bounded by the maximum of
\hyperref[eq:3.7]{(3.7)}, divided by 240. The finite checks bound that maximum by 596 after
the certified exclusions. Consequently \[
\operatorname{out}(h)\le\frac{596}{240}=\frac{149}{60},
\qquad\operatorname{ch}(h)\ge\frac{25}{2}+\frac1{60}.
\tag{3.8}\label{eq:3.8}
\]

A root exclusion identifies the common neighbors of an order-nine link
and \(L\) along their center edge. Verified root automorphisms and
remaining permutations cover all identifications. A case is excluded by
a forced violation of \(\mathcal C\), or a certificate satisfying \hyperref[eq:2.1]{(2.1)}
for the centers' spokes. Both centers' complete induced links are
specified. Every triangle using a deleted spoke is therefore determined
by the corresponding complete link. The forced-edge and
triangle-coverage conditions of \hyperlink{statement-lemma-A-1}{Lemma~A.1} give validity in every
permitted containing graph; unknown private cross edges cannot create
another triangle on a center's spoke. Retained roots stay allowed; they
are not assumed globally realizable. \hyperref[appendix-a.-finite-domains-and-local-graph-data]{Appendix~A} and \hyperref[verification-scope-and-a-remaining-question]{Section~6} specify
the verification boundary.

\Needspace{5\baselineskip}
\subsection{3.4. Proof of Theorem 1.1}\label{proof-of-theorem-1.1}

Suppose an input satisfying the hypotheses has \(\Phi(G)>0\). Apply
reductions while possible and discard isolates. Lemmas 2.2 and 2.4 yield
a nonempty reduction-free \(J\in\mathcal C\) with positive \(\Phi\) and
\(\operatorname{mad}(J)\le25/2\). Every vertex finishes with charge at
least \(25/2\), and every high vertex strictly more. Conservation and
the density bound exclude any high vertex. \hyperlink{statement-lemma-2-5}{Lemma~2.5} then excludes
degrees nine, ten and eleven. All remaining degrees are twelve or
thirteen. \hyperref[appendix-b.-the-degree-twelvethirteen-endpoint]{Appendix~B} excludes a nonempty reduced graph of this form with
average degree at most \(25/2\). This proves \hyperlink{statement-theorem-1-1}{Theorem~1.1}, including the
endpoint. \hfill\(\square\)

\Needspace{5\baselineskip}
\section{4. Localization to the Hifumi
core}\label{localization-to-the-hifumi-core}

This section concerns the larger class \(\mathcal D\) and does not use
the density theorem.

Recall the distinction between \(K\) and \(G[K]\) illustrated in
\hyperref[fig-hifumi-core]{Figure 1}; the reduction below constructs a
possibly proper subgraph of \(G[K]\).

\par\addvspace{7pt}\Needspace{5\baselineskip}\noindent

\textbf{\statementtarget{statement-lemma-4-1}Lemma 4.1 (links of arbitrary order).} Every finite graph \(L\)
of minimum degree at least two and maximum degree at most four, with at
most two vertices of degree three and at most two of degree four, is
WKE. It may be disconnected.

\hyperref[appendix-c.-wke-witnesses-for-links-of-arbitrary-order]{Appendix~C} proves this using finite base graphs and a
certificate-preserving two-vertex subdivision. The order of \(L\) is not
bounded.

\par\addvspace{7pt}\Needspace{5\baselineskip}\noindent

\Needspace{11\baselineskip}
\textbf{\statementtarget{statement-lemma-4-2}Lemma 4.2 (residual structure).} In every nonempty
reduction-free \(G\in\mathcal D\), with isolated vertices removed, \[
d_{H_3(G)}(v)\ge4,\qquad
d_{H_3(G)}(v)\equiv0\pmod2
\quad(v\in V(G)).
\tag{4.1}\label{eq:4.1}
\] \emph{Proof.} \hyperlink{statement-lemma-2-4}{Lemma~2.4} restricts degrees in \(L_v\) to two, three or
four, with at most two fours. The number of threes equals
\(d_{H_3(G)}(v)\) and is even by the handshaking lemma. If it were zero
or two, \hyperlink{statement-lemma-4-1}{Lemma~4.1} and \hyperlink{statement-lemma-2-3}{Lemma~2.3} would reduce the spokes at \(v\). \hfill\(\square\)

\textbf{Proof of Theorem 1.2.} Fix a peeling order of
\(V(G)\setminus K\) in the initial \(H_{\ge3}(G)\). At every stage the
current \(J\) is a subgraph of \(G\) in \(\mathcal D\), and \[
H_{\ge3}(J)\subseteq H_{\ge3}(G)[V(J)].
\tag{4.2}\label{eq:4.2}
\] Repeatedly remove codegree-zero and codegree-one edges, and discard
isolates. If a vertex outside \(K\) remains, choose the earliest in the
fixed peeling order, say \(v\). It has at most three current
\(H_{\ge3}(J)\)-neighbors. All incident codegrees are now two, three or
four. The number of threes is even and at most three, hence at most two;
the number of fours is at most two. \hyperlink{statement-lemma-4-1}{Lemma~4.1} supplies a WKE witness in
\(L_v\), and \hyperlink{statement-lemma-2-3}{Lemma~2.3} removes its spokes. Remove the now isolated
center and continue. Each nontrivial step removes edges, so the process
terminates with \(R\) supported on \(K\). Equation \hyperref[eq:2.2]{(2.2)} telescopes to
\hyperref[eq:1.3]{(1.3)}. \hfill\(\square\)

Equation \hyperref[eq:4.2]{(4.2)} would fail with \(H_3\) in place of \(H_{\ge3}\):
deletion can reduce a codegree from four to three. This explains the use
of \emph{at least} three triangles. Further, Tuza's inequality for
\(G[K]\) alone does not imply it for the survivor \(R\). \hyperlink{statement-corollary-1-3}{Corollary~1.3}
uses the hereditary sufficient hypotheses of Theorem 1.1.

For \(G\in\mathcal C\), the quantitative gap also implies that if
\(\Phi(G)>0\), then \(K\ne\varnothing\) and \[
\operatorname{mad}(G[K])\ge\frac{25}{2}+\eta.
\tag{4.3}\label{eq:4.3}
\] An edge-minimal counterexample subgraph of the survivor is
reduction-free and satisfies both \hyperref[eq:4.1]{(4.1)} and Theorem 1.4. These
statements concern a reduced subgraph, not the ordinary average degree
of every input counterexample.

\Needspace{5\baselineskip}
\Needspace{23\baselineskip}
\section{5. A uniform positive
surplus}\label{a-uniform-positive-surplus}

Let \(G\) be a nonempty reduction-free graph in \(\mathcal C\) without
isolates, with no density hypothesis. Apply \hyperref[discharging-and-the-density-bound]{Section~3} and put \[
s(v)=\operatorname{ch}(v)-25/2,\qquad S=\sum_{v\in V(G)}s(v).
\tag{5.1}\label{eq:5.1}
\] All surpluses are nonnegative. High vertices have explicit positive
lower bounds, but need not themselves form a positive fraction of the
graph. \hyperref[appendix-d.-the-quantitative-surplus-count]{Appendix~D} controls degree-twelve and degree-thirteen vertices
through the local profiles of \hyperref[appendix-b.-the-degree-twelvethirteen-endpoint]{Appendix~B}, and proves \[
|V(G)|\le\frac{56775}{2}S.
\tag{5.2}\label{eq:5.2}
\] Conservation yields \[
\frac{2|E(G)|}{|V(G)|}
=\frac{25}{2}+\frac{S}{|V(G)|}
\ge\frac{25}{2}+\frac{2}{56775}.
\tag{5.3}\label{eq:5.3}
\] An input in \(\mathcal C\) below this maximum-average-degree
threshold that violated Tuza's inequality would reduce to a nonempty
subgraph contradicting \hyperref[eq:5.3]{(5.3)}. Applying the same argument to
\(R\subseteq G[K]\) gives the localized version. This proves Theorem
1.4. The strict endpoint in \hyperref[eq:1.4]{(1.4)} must not be replaced by a closed
endpoint at \(\eta\). \hfill\(\square\)

\Needspace{5\baselineskip}
\clearpage
\section{6. Verification, scope and a remaining
question}\label{verification-scope-and-a-remaining-question}

\Needspace{5\baselineskip}
\subsection{6.1. Finite premises and analytic
implications}\label{finite-premises-and-analytic-implications}

\hyperref[tab-verification]{Table 3} records the finite obligations and
the analytic arguments that use them.

\Needspace{10\baselineskip}
\begin{papertable}
\caption{Finite verification obligations and the analytic arguments needed beyond
the checks.}
\label{tab-verification}
\hypertarget{tab-verification}{}
\vspace{4pt}\centering\small\renewcommand{\arraystretch}{1.22}
\begin{tabular}{@{}>{\raggedright\arraybackslash}p{19.630mm}@{\hspace{3mm}}>{\raggedright\arraybackslash}p{37.750mm}@{\hspace{3mm}}>{\raggedright\arraybackslash}p{36.240mm}@{\hspace{3mm}}>{\raggedright\arraybackslash}p{57.380mm}@{}}
\toprule
\textbf{ID} & \textbf{Finite obligation} & \textbf{Mathematical use} & \textbf{Argument outside the finite check} \\
\midrule
F001 & Small links and their certificates & \hyperlink{statement-lemma-2-5}{Lemma~2.5}(a); \hyperref[appendix-e.-an-ordinary-core-sufficient-condition]{Appendix~E} & WKE/cone reductions in arbitrary containing graphs \\
F002--F004 & Orders 9--14, root maps, shared-center conditions & \hyperlink{statement-lemma-2-5}{Lemma~2.5}(b)--(f); Appendices B and D & Actual neighborhoods give tested identifications; \hyperlink{statement-lemma-A-1}{Lemma~A.1} and Appendix
A.2 justify ambient validity \\
F005 & Degree-fourteen and degree-sixteen integer bounds & \hyperref[discharging-and-the-density-bound]{Section~3} & Conservation and the estimate for unbounded degree \\
F006 & Order-fifteen domain, exclusions, interfaces and \hyperref[eq:3.7]{(3.7)} & Equation \hyperref[eq:3.8]{(3.8)} & The actual paid set is contained in the upper-bound domain \\
F007--F008 & WKE base graphs and complete parity domains & \hyperlink{statement-lemma-4-1}{Lemma~4.1} & Subdivision and disjoint union, \hyperref[appendix-c.-wke-witnesses-for-links-of-arbitrary-order]{Appendix~C} \\
\bottomrule
\end{tabular}
\end{papertable}

The finite checks establish literal witnesses and domain coverage. The
proof is not a claim about many sampled graphs: completeness and the
actual-to-local mapping are indispensable. A retained configuration need
not extend to a graph, and a failed search is never used as a
certificate of nonexistence.

The order-fifteen calculation uses two different independent-set
maximizers, but they share parsing, interface application and high-set
enumeration. Their agreement is not two fully independent proofs.
Separate interface checks and the actual-to-local argument address that
boundary. Graph generation uses nauty \hyperlink{bib-6}{[6]}; its routines, the
monotone pruning conditions, compiler and integer arithmetic
implementation remain part of the computational trust boundary. This is
not a proof-assistant formalization.

\textbf{Verification materials.} The companion package contains the fixed
inputs and source checkers for F001--F008. Its \texttt{START\_HERE.html}
guide maps each check to the paper, gives environment and resource
instructions, and explains how to run
\mbox{\texttt{python3 -B verify.py -{}-full}}. Each run builds fresh tools and
records its results and file identities. The finite checks do not replace
the analytic arguments above. The verification package has been deposited as an
unpublished draft on Zenodo and will be made publicly available after
this article is announced on arXiv. Reserved DOI:
\href{https://doi.org/10.5281/zenodo.24047493}{\nolinkurl{10.5281/zenodo.24047493}}.

\Needspace{5\baselineskip}
\Needspace{18\baselineskip}
\subsection{6.2. What a counterexample would still
require}\label{what-a-counterexample-would-still-require}

\par\addvspace{7pt}\Needspace{5\baselineskip}\noindent

\textbf{\statementtarget{statement-corollary-6-1}Corollary 6.1.} If Tuza's conjecture is false, there is an
edge-minimal counterexample \(J\), without isolates, for which at least
one of the following holds:

\begin{enumerate}
\def\labelenumi{\arabic{enumi}.}
\tightlist
\item
  \(b(J)\ge5\);
\item
  \(b(J)\le4\), and some component of \(H_4(J)\) has at least three
  edges;
\item
  \(J\in\mathcal C\) and \(2|E(J)|/|V(J)|\ge25/2+\eta\).
\end{enumerate}

In the third case every vertex also satisfies \hyperref[eq:4.1]{(4.1)}.

\emph{Proof.} If the first two alternatives fail, \(J\in\mathcal C\).
Edge minimality and \hyperlink{statement-lemma-2-2}{Lemma~2.2} exclude every reduction. \hyperlink{statement-theorem-1-4}{Theorem~1.4} and
\hyperlink{statement-lemma-4-2}{Lemma~4.2} apply. \hfill\(\square\)

\Needspace{11\baselineskip}
Even with isolated vertices excluded, this ordinary-average-degree
conclusion cannot be transferred to arbitrary counterexamples. If \(J\)
is a counterexample without isolates, adjoining \(m\) disjoint copies of
\(K_2\) changes neither triangle parameter and preserves membership in
\(\mathcal C\), when applicable, while \[
\frac{2|E(J\sqcup mK_2)|}{|V(J\sqcup mK_2)|}
=\frac{2|E(J)|+2m}{|V(J)|+2m}\longrightarrow1.
\] These enlarged graphs have no isolated vertices but are not
edge-minimal counterexamples.

\par\addvspace{7pt}\Needspace{5\baselineskip}\noindent

\textbf{\statementtarget{statement-question-6-2}Question 6.2.} Does every graph \(G\in\mathcal C\) satisfy
\(\tau(G)\le2\nu(G)\), without a density condition?

A positive answer would eliminate the third alternative of Corollary
6.1, but would not settle the first two or the general conjecture. We do
not resolve \hyperlink{statement-question-6-2}{Question~6.2}, classify all reduced graphs, or prove that the
thresholds are optimal. The question identifies a specific extension,
not a claim that only one obstacle remains to all approaches.

\Needspace{5\baselineskip}
\section{Appendix A. Finite domains and local graph
data}\label{appendix-a.-finite-domains-and-local-graph-data}

\Needspace{5\baselineskip}
\subsection{A.1. Small components and the low-order
domain}\label{a.1.-small-components-and-the-low-order-domain}

The small-component domain comprises 1,473 connected graphs of order at
most eight, maximum degree at most four and at most two vertices of
degree four. There are 1,471 WKE witnesses and two direct cone
witnesses, checked against \hyperref[eq:2.1]{(2.1)} or \hyperref[eq:2.3]{(2.3)}. Every triangle using a spoke
of one link component stays in that component's cone; exterior edges
therefore cannot invalidate the reduction. This proves the premise used
in \hyperlink{statement-lemma-2-5}{Lemma~2.5}(a).

The domain sizes for connected links of minimum degree two, maximum
degree four and at most two degree-four vertices are listed in
\hyperref[tab-link-counts]{Table 4}.

\Needspace{10\baselineskip}
\begin{papertable}
\caption{Connected link domains of orders nine through fourteen under the stated
degree constraints.}
\label{tab-link-counts}
\hypertarget{tab-link-counts}{}
\vspace{4pt}\centering\small\renewcommand{\arraystretch}{1.22}
\begin{tabular}{@{}>{\raggedright\arraybackslash}p{17.040mm}@{\hspace{3mm}}>{\raggedright\arraybackslash}p{19.880mm}@{\hspace{3mm}}>{\raggedright\arraybackslash}p{19.880mm}@{\hspace{3mm}}>{\raggedright\arraybackslash}p{19.880mm}@{\hspace{3mm}}>{\raggedright\arraybackslash}p{19.880mm}@{\hspace{3mm}}>{\raggedright\arraybackslash}p{21.300mm}@{\hspace{3mm}}>{\raggedright\arraybackslash}p{24.140mm}@{}}
\toprule
\textbf{Order} & \textbf{9} & \textbf{10} & \textbf{11} & \textbf{12} & \textbf{13} & \textbf{14} \\
\midrule
Graphs & 1,691 & 7,612 & 34,721 & 164,523 & 791,016 & 3,885,906 \\
\bottomrule
\end{tabular}
\end{papertable}

The first four orders total 208,547. Of these, 207,696 have WKE
witnesses and 720 have cones already violating \(\mathcal C\). Of the
remaining 131 types, cone and one-sided contact certificates exclude 106
and 8 respectively, leaving seventeen templates. The one-sided
exclusions cover all 183 relevant contact patterns.

These counts specify finite obligations. Merely counting saved records
does not establish completeness: the graph domain must be regenerated
with the stated constraints and every exclusion checked.

\Needspace{5\baselineskip}
\subsection{A.2. Adjacent centers and arbitrary exterior
graphs}\label{a.2.-adjacent-centers-and-arbitrary-exterior-graphs}

For adjacent centers \(u,v\), common neighbors are identified under
every allowed bijection, modulo verified root automorphisms. Private
neighborhoods are disjoint by definition. Adjacencies determined by the
induced links must agree. Other edges and exterior vertices are not
assumed absent. A reduction must cover all triangles using its deleted
spokes in every permitted completion.

\par\addvspace{7pt}\Needspace{5\baselineskip}\noindent

\Needspace{12\baselineskip}
\textbf{\statementtarget{statement-lemma-A-1}Lemma A.1 (ambient validity).} Suppose a local record is
identified injectively with vertices of a containing graph. Let
\(E_{\rm rec}\) be its forced-present edges, and let
\(F\subseteq E_{\rm rec}\) be nonempty. Suppose \(\mathcal S\) is an
edge-disjoint triangle family using \(E_{\rm rec}\),
\(X\subseteq E_{\rm rec}\), \(E(\mathcal S)\setminus F\subseteq X\), and
\(|X|\le2|\mathcal S|\). For every \(xy\in F\setminus X\), require that
the record accounts for all possible common neighbors of \(x,y\), and
all corresponding triangles in every permitted completion, and that
\(X\) meets each such triangle. Then \((F,\mathcal S,X)\) satisfies
\hyperref[eq:2.1]{(2.1)} in every permitted completion.

\emph{Proof.} The packing and cover use forced edges, so their
existence, disjointness and size inequalities persist. Every triangle
using \(F\setminus X\) is among those tested; a triangle using
\(F\cap X\) is met by that edge itself. These are exactly the remaining
conditions of \hyperref[eq:2.1]{(2.1)}. \hfill\(\square\)

When both centers' complete induced links are specified, every triangle
through a center's spoke is known. These are the deletion edges in the
two-center certificates. Edges between their private neighborhoods and
arbitrary vertices outside both neighborhoods cannot create an
additional such triangle. Literal checks verify the forced-edge, packing
and cover conditions of Lemma A.1.

A cone containing an edge in five triangles cannot occur when \(b\le4\).
Nor can a component of codegree-four edges have more than two edges in
\(\mathcal C\): the four triangles witnessing each edge persist in a
containing graph. These are forced violations, not failed searches.

One-sided contacts do not bound the opposite vertex degree. Let
\(C=N(u)\cap N(v)\) at a codegree-four edge \(uv\). Its four vertices
are fixed. If \(u\) is the central endpoint of a \(P_3\) in \(H_4\),
then \(v\) is a leaf there and all other edges at \(v\) have codegree at
most three. For \(c\in C\), the codegree of \(vc\) already counts \(u\)
and its \(d_{G[C]}(c)\) neighbors in \(C\); at most \(2-d_{G[C]}(c)\)
additional common neighbors are possible. If \(u\) is a leaf, \(v\) may
have one other codegree-four edge, allowing at most one extra contact
unit. Negative capacities force exclusion. Every private vertex meeting
\(C\) consumes a positive amount of this bounded capacity. Enumerating
multisets of their nonempty contact subsets thus gives a finite domain.
Private vertices not meeting \(C\), and unknown edges, are still allowed
in the certificate model; their total number is not bounded by this
enumeration.

\Needspace{10\baselineskip}
For clarity, the one-sided records use a more restricted deletion set.
Write \[
P=N(u)\setminus(C\cup\{v\}),\qquad
Z_+=\{z\in N(v)\setminus(C\cup\{u\}):N(z)\cap C\ne\varnothing\},
\] and let \(Z_0\) be the remaining private neighbors of \(v\). The
complete induced link at \(u\) is recorded. Each vertex of \(Z_+\) is
recorded separately with its exact nonempty contact set in \(C\);
repeated contact sets represent distinct vertices. The capacity
enumeration includes every such multiset, without requiring that
capacities be exhausted. The number of vertices of \(Z_0\) is
unrestricted.

Every one-sided certificate uses exactly \[
F=\{ux:x\in N(u)\}\ \cup\ \{vc:c\in C\}.
\] In particular, no spoke from \(v\) to \(Z_+\) or \(Z_0\) belongs to
\(F\). The reduced graph deletes \(F\cup X\), as in Lemma 2.2. The
triangles on the first group of edges are fixed by the complete induced
link at \(u\).

\Needspace{13\baselineskip}
For \(c\in C\), the entire common-neighbor set of \(v,c\)
is \[
N(v)\cap N(c)=\{u\}\ \cup\ (C\cap N(c))
\cup\{z\in Z_+:c\in N(z)\}.
\] Every term is recorded: a vertex of \(P\) is not adjacent to \(v\),
one of \(Z_0\) is not adjacent to \(c\), and a vertex outside both
neighborhoods is not adjacent to \(v\). Thus all triangles on these
deletion edges are also fixed. Unspecified edges cannot add another
triangle using \(F\). The acceptance predicate checks this shape of
\(F\), the complete link and exact contact data, and all conditions of
\hyperlink{statement-lemma-A-1}{Lemma~A.1} using forced-present edges. This is the all-completions
argument; the capacity bound alone is not used as a lift.

Among the seventeen low-order templates there are 1,342 compatible
root-orbit/permutation cases. Literal reductions exclude 1,182 and class
violations exclude 139. All 21 retained cases remain possible in later
tests. No global extension is assumed.

\Needspace{5\baselineskip}
\subsection{A.3. Shared-center consistency and degree
twelve}\label{a.3.-shared-center-consistency-and-degree-twelve}

All but one retained degree-twelve template already have at most two
possible neighbors of degree at most twelve. The exception \(T_C\) has
candidate labels \(0,5,9\). Labels 5 and 9 are not adjacent; 0 is
adjacent to both. The pair table under all retained maps and
automorphisms forbids 0 and 5 being low simultaneously, and likewise 0
and 9.

The test retains third-center identity. If \(v,w\) are the neighbors at
labels \(i,j\) in the link of \(u\), the \(u,v\) correspondence locates
\(w\) in the link of \(v\), and the \(u,w\) correspondence locates \(v\)
in the link of \(w\). A retained \(v,w\) correspondence must identify
the same vertex \(u\). None does in the two forbidden cases. All 45
directed root-pair keys and 88 common-neighbor maps are tested. Other
globally inconsistent cases may remain, safely enlarging the
necessary-condition domain. Thus each degree-twelve vertex has at most
two neighbors of degree at most twelve.

Every eleven-element subset of every retained link of a degree-thirteen
vertex contains a vertex with at least three neighbors in that subset.
If eleven neighbors of a degree-thirteen vertex all had degree twelve,
this would violate their preceding bound. Hence a vertex of degree
thirteen has at most ten degree-twelve neighbors.

The root tables also establish \hyperlink{statement-lemma-2-5}{Lemma~2.5}(b)--(d) and (f). A designation
absent from every retained pair of degree at most thirteen must be high.
The degree-fourteen table rules out contacts with degree-nine vertices
designating the center, and allows contacts with degree-ten vertices
only between \(T_R\) at label 4 (or its image 9) and \(J_{14}\) at label
1 (or its image 6). Pairs whose centers both have degree nine or ten,
including the undesignated label of \(N_4\), are excluded. The nine
designations at a degree-eleven vertex are selected from its eligible
labels, each with the two-high property.

\Needspace{5\baselineskip}
\clearpage
\subsection{A.4. Explicit graphs used in the analytic
arguments}\label{a.4.-explicit-graphs-used-in-the-analytic-arguments}

These graph names are labels, not additional hypotheses. Every graph has
vertices \mbox{\(0,1,\ldots,n-1\)}; the unordered edge pairs below fix its
labeling. The full collection and root maps are separate finite inputs.

\noindent\begin{minipage}{\linewidth}
\noindent\begin{minipage}[c]{.16\linewidth}
\textbf{\(N_1\), \(n=9\).}
\end{minipage}\hfill\begin{minipage}[c]{.82\linewidth}
\setlength{\abovedisplayskip}{2pt}\setlength{\belowdisplayskip}{2pt}\setlength{\abovedisplayshortskip}{2pt}\setlength{\belowdisplayshortskip}{2pt}
\[
\begin{gathered}
(0,4),(0,6),(0,8),(1,4),(1,6),(1,8),(2,5),\\
(2,7),(2,8),(3,5),(3,7),(3,8),(4,6),(5,7)
\end{gathered}
\]
\end{minipage}\par\vspace{5pt}

\noindent\begin{minipage}[c]{.16\linewidth}
\textbf{\(N_2\), \(n=9\).}
\end{minipage}\hfill\begin{minipage}[c]{.82\linewidth}
\setlength{\abovedisplayskip}{2pt}\setlength{\belowdisplayskip}{2pt}\setlength{\abovedisplayshortskip}{2pt}\setlength{\belowdisplayshortskip}{2pt}
\[
\begin{gathered}
(0,4),(0,7),(0,8),(1,4),(1,7),(1,8),(2,5),\\
(2,6),(2,7),(3,5),(3,6),(3,8),(4,6),(5,8)
\end{gathered}
\]
\end{minipage}\par\vspace{5pt}

\noindent\begin{minipage}[c]{.16\linewidth}
\textbf{\(N_3\), \(n=9\).}
\end{minipage}\hfill\begin{minipage}[c]{.82\linewidth}
\setlength{\abovedisplayskip}{2pt}\setlength{\belowdisplayskip}{2pt}\setlength{\abovedisplayshortskip}{2pt}\setlength{\belowdisplayshortskip}{2pt}
\[
\begin{gathered}
(0,4),(0,7),(0,8),(1,4),(1,7),(1,8),(2,5),\\
(2,6),(2,8),(3,5),(3,6),(3,8),(4,6),(5,7)
\end{gathered}
\]
\end{minipage}\par\vspace{5pt}

\noindent\begin{minipage}[c]{.16\linewidth}
\textbf{\(N_4\), \(n=9\).}
\end{minipage}\hfill\begin{minipage}[c]{.82\linewidth}
\setlength{\abovedisplayskip}{2pt}\setlength{\belowdisplayskip}{2pt}\setlength{\abovedisplayshortskip}{2pt}\setlength{\belowdisplayshortskip}{2pt}
\[
\begin{gathered}
(0,4),(0,5),(0,7),(1,4),(1,5),(1,7),(2,6),\\
(2,7),(2,8),(3,6),(3,7),(3,8),(4,8),(5,8)
\end{gathered}
\]
\end{minipage}\par\vspace{5pt}

\noindent\begin{minipage}[c]{.16\linewidth}
\textbf{\(T_P\), \(n=10\).}
\end{minipage}\hfill\begin{minipage}[c]{.82\linewidth}
\setlength{\abovedisplayskip}{2pt}\setlength{\belowdisplayskip}{2pt}\setlength{\abovedisplayshortskip}{2pt}\setlength{\belowdisplayshortskip}{2pt}
\[
\begin{gathered}
(0,7),(0,8),(0,9),(1,5),(1,6),(1,9),(2,4),\\
(2,6),(2,8),(3,4),(3,5),(3,7),(4,9),(5,8),\\
(6,7)
\end{gathered}
\]
\end{minipage}\par\vspace{5pt}

\noindent\begin{minipage}[c]{.16\linewidth}
\textbf{\(T_R\), \(n=10\).}
\end{minipage}\hfill\begin{minipage}[c]{.82\linewidth}
\setlength{\abovedisplayskip}{2pt}\setlength{\belowdisplayskip}{2pt}\setlength{\abovedisplayshortskip}{2pt}\setlength{\belowdisplayshortskip}{2pt}
\[
\begin{gathered}
(0,2),(0,3),(0,4),(1,2),(1,3),(1,4),(2,3),\\
(4,9),(5,7),(5,8),(5,9),(6,7),(6,8),(6,9),\\
(7,8)
\end{gathered}
\]
\end{minipage}\par\vspace{5pt}

\noindent\begin{minipage}[c]{.16\linewidth}
\textbf{\(T_A\), \(n=12\).}
\end{minipage}\hfill\begin{minipage}[c]{.82\linewidth}
\setlength{\abovedisplayskip}{2pt}\setlength{\belowdisplayskip}{2pt}\setlength{\abovedisplayshortskip}{2pt}\setlength{\belowdisplayshortskip}{2pt}
\[
\begin{gathered}
(0,5),(0,10),(0,11),(1,6),(1,8),(1,9),(2,6),\\
(2,8),(2,9),(3,7),(3,9),(4,7),(4,10),(4,11),\\
(5,10),(5,11),(6,8)
\end{gathered}
\]
\end{minipage}\par\vspace{5pt}

\noindent\begin{minipage}[c]{.16\linewidth}
\textbf{\(T_C\), \(n=12\).}
\end{minipage}\hfill\begin{minipage}[c]{.82\linewidth}
\setlength{\abovedisplayskip}{2pt}\setlength{\belowdisplayskip}{2pt}\setlength{\abovedisplayshortskip}{2pt}\setlength{\belowdisplayshortskip}{2pt}
\[
\begin{gathered}
(0,5),(0,8),(0,9),(1,6),(1,7),(1,8),(2,6),\\
(2,7),(2,8),(3,9),(3,10),(3,11),(4,9),(4,10),\\
(4,11),(5,10),(5,11),(6,7)
\end{gathered}
\]
\end{minipage}\par\vspace{5pt}

\noindent\begin{minipage}[c]{.16\linewidth}
\textbf{\(T_D\), \(n=12\).}
\end{minipage}\hfill\begin{minipage}[c]{.82\linewidth}
\setlength{\abovedisplayskip}{2pt}\setlength{\belowdisplayskip}{2pt}\setlength{\abovedisplayshortskip}{2pt}\setlength{\belowdisplayshortskip}{2pt}
\[
\begin{gathered}
(0,5),(0,9),(0,11),(1,6),(1,8),(1,9),(2,6),\\
(2,8),(2,9),(3,7),(3,10),(3,11),(4,7),(4,10),\\
(4,11),(5,10),(6,8),(7,11)
\end{gathered}
\]
\end{minipage}\par\vspace{5pt}

\noindent\begin{minipage}[c]{.16\linewidth}
\textbf{\(J_{14}\), \(n=14\).}
\end{minipage}\hfill\begin{minipage}[c]{.82\linewidth}
\setlength{\abovedisplayskip}{2pt}\setlength{\belowdisplayskip}{2pt}\setlength{\abovedisplayshortskip}{2pt}\setlength{\belowdisplayshortskip}{2pt}
\[
\begin{gathered}
(0,6),(0,7),(1,7),(1,12),(1,13),(2,8),(2,10),\\
(2,12),(3,8),(3,10),(3,12),(4,9),(4,11),(4,13),\\
(5,9),(5,11),(5,13),(6,12),(6,13),(8,10),(9,11)
\end{gathered}
\]
\end{minipage}\par\vspace{5pt}
\end{minipage}\par

\Needspace{8\baselineskip}
The three degree-twelve link types with two degree-twelve neighbors are
\(T_A,T_C,T_D\). For \(T_C\) the two labels must be 5 and 9. Each
corresponding neighbor has type \(T_C\) or \(T_D\) with the original
center at its label 0; the retained shared-center profiles then force it
to have at most one neighbor of degree at most twelve. In \(T_D\), label
11 has link degree four; the low pair table excludes degrees at most
twelve there and positive certificates exclude degree thirteen. In
\(T_A\), the two degree-twelve neighbor labels are 3 and 7, and every
other label is adjacent to at most one of them. Finally every retained
link of a degree-thirteen vertex has at least eighteen edges. These
finite facts are used at equality and, without a global equality
assumption, in \hyperref[appendix-d.-the-quantitative-surplus-count]{Appendix~D}.

\Needspace{5\baselineskip}
\Needspace{28\baselineskip}
\section{Appendix B. The degree-twelve/thirteen
endpoint}\label{appendix-b.-the-degree-twelvethirteen-endpoint}

\par\addvspace{7pt}\Needspace{5\baselineskip}\noindent

\textbf{\statementtarget{statement-lemma-B-1}Lemma B.1.} No nonempty reduction-free graph in \(\mathcal C\),
all of whose degrees are twelve or thirteen, has average degree at most
\(25/2\).

\emph{Proof.} Let \(n_{12},n_{13}\) count the two degree classes. Lemma
2.5(e) gives \[
10n_{12}\le e(V_{12},V_{13})\le10n_{13}.
\tag{B.1}\label{eq:B.1}
\] The average-degree condition gives \(n_{13}\le n_{12}\). Thus the
counts are equal and positive; every degree-twelve vertex has exactly
two degree-twelve neighbors, and every degree-thirteen vertex has
exactly ten degree-twelve neighbors.

Appendix A.4 leaves \(T_A,T_C,T_D\) for the degree-twelve vertices. A
\(T_C\) vertex forces its two degree-twelve neighbors to have at most
one low neighbor, contradicting equality. A \(T_D\) vertex has its
label-11 neighbor of degree at least fourteen, also impossible.
Therefore every degree-twelve vertex has type \(T_A\), with its two
degree-twelve neighbors at labels 3 and 7.

Fix a degree-thirteen center \(h\). Its ten degree-twelve neighbors
induce maximum degree at most one in \(L_h\): for each such neighbor
\(v\), the label of \(h\) in its \(T_A\)-link is outside \(\{3,7\}\) and
is adjacent to at most one of them. Thus these ten vertices induce at
most five edges. Yet \(L_h\) has at least eighteen edges, and its
remaining three vertices have degree sum at most eleven, since at most
two link vertices have degree four. Removing them leaves at least seven
edges on the ten degree-twelve vertices, a contradiction. This proof
does not use Theorem 1.1. \hfill\(\square\)

\Needspace{5\baselineskip}
\Needspace{14\baselineskip}
\section{Appendix C. WKE witnesses for links of arbitrary
order}\label{appendix-c.-wke-witnesses-for-links-of-arbitrary-order}

\Needspace{5\baselineskip}
\subsection{C.1. Subdivision preserves a
witness}\label{c.1.-subdivision-preserves-a-witness}

Replace an edge \(uv\) of \(L\) by \(ua,ab,bv\), with \(a,b\) new. Given
a WKE witness \((M,Q)\):

\begin{itemize}
\tightlist
\item
  if \(uv\in M\), replace it in \(M\) by \(ua,bv\) and add \(a\) to
  \(Q\);
\item
  if \(uv\notin M\), at least one endpoint is in \(Q\). Add \(ab\) to
  \(M\); add \(b\) to \(Q\) if \(u\in Q\), and otherwise add \(a\).
\end{itemize}

In both cases the new set is a matching, the new \(Q\) covers every edge
outside it, and both sizes increase by one. All other incidences are
unchanged. Repetition lifts a witness through any even number of added
internal vertices. Disjoint unions combine witnesses additively.

\hyperref[fig-wke-subdivision]{Figure 3} illustrates the two witness
updates. The displayed triangles are examples; the preceding
verification applies to an arbitrary graph carrying a WKE witness.

\begin{paperfigure}
\centering
\includegraphics[width=\linewidth]{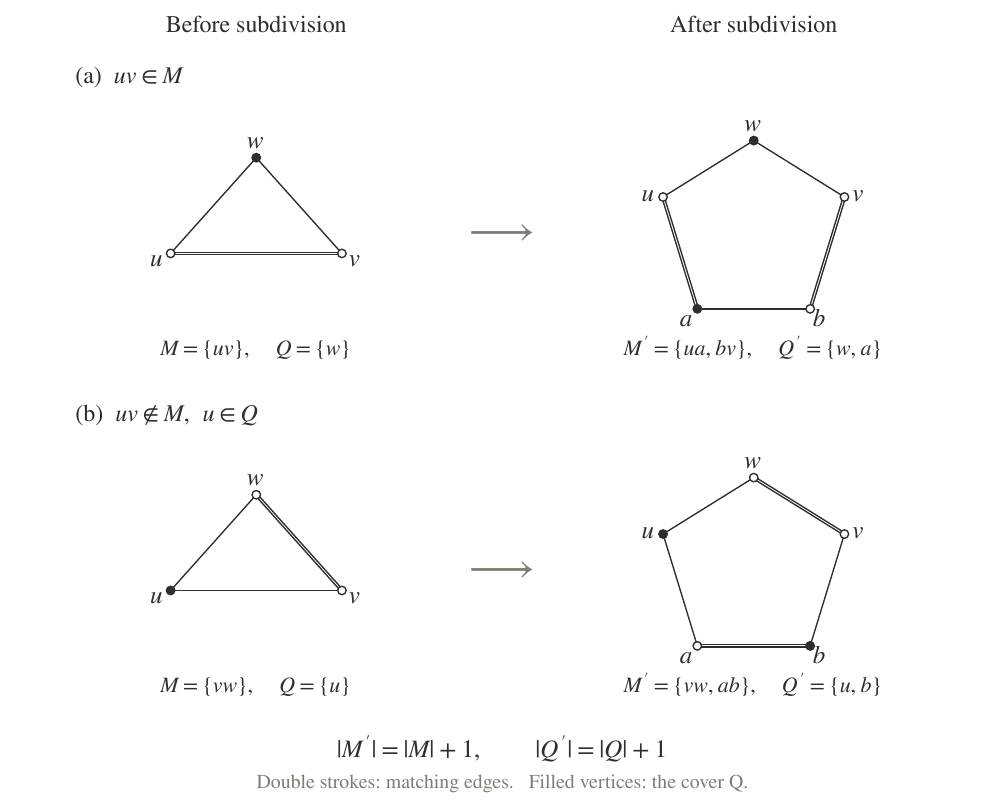}
\caption{Two examples of the witness-preserving subdivision in Appendix C.1.
Replace \(uv\) by \(ua,ab,bv\). Double strokes mark single matching
edges; filled vertices form the cover \(Q\) of the edges outside \(M\),
or the corresponding cover \(Q'\) after subdivision. In (a),
\(uv\in M\): replace it by \(ua,bv\) and add \(a\) to \(Q\). In (b),
\(uv\notin M\) and \(u\in Q\): add \(ab\) to \(M\) and \(b\) to \(Q\).
The case \(u\notin Q\), \(v\in Q\) is obtained by exchanging the
endpoints. In both examples the matching and cover sizes each increase
by one. The general lifting argument, including repeated subdivisions,
is given in the text.}
\label{fig-wke-subdivision}
\hypertarget{fig-wke-subdivision}{}
\end{paperfigure}

\Needspace{5\baselineskip}
\Needspace{9\baselineskip}
\subsection{C.2. At most two branch
vertices}\label{c.2.-at-most-two-branch-vertices}

Call a vertex of degree greater than two a branch vertex in this
appendix. A connected graph with no branch is a cycle, reduced by the
preceding operation to \(C_3\) or \(C_4\). Suppressing maximal paths
with internal degree-two vertices gives a multigraph on the branches.
Loops contribute twice to a branch degree.

One degree-four branch gives two loops. For two degree-four branches,
with \(p>0\) connecting paths, \(2l_i+p=4\) gives either four connecting
paths, or two paths and a loop at each branch. For two degree-three
branches, \(2l_i+p=3\) gives three connecting paths, or one path and a
loop at each branch.

Retain a loop of length three when odd and four when even. Retain a
nonloop path of length one if it is an original direct edge, two if
even, and three if odd but not direct. Simplicity permits at most one
direct edge between a pair of branches. The complete finite collections
are listed in \hyperref[tab-two-branch]{Table 5}.

\Needspace{10\baselineskip}
\begin{papertable}
\caption{Complete labeled parity/path domains with at most two branch vertices.}
\label{tab-two-branch}
\hypertarget{tab-two-branch}{}
\vspace{4pt}\centering\small\renewcommand{\arraystretch}{1.22}
\begin{tabular}{@{}>{\raggedright\arraybackslash}p{117.750mm}@{\hspace{3mm}}>{\raggedright\arraybackslash}p{39.250mm}@{}}
\toprule
\textbf{Topology} & \textbf{Labeled parity/path records} \\
\midrule
Cycles & 2 \\
One degree-four branch & 4 \\
Four paths between degree-four branches & 48 \\
Two connecting paths and two loops, degree four & 32 \\
Three paths between degree-three branches & 20 \\
One connecting path and two loops, degree three & 12 \\
\textbf{Total} & \textbf{118} \\
\bottomrule
\end{tabular}
\end{papertable}

Each base graph has a literal WKE witness. These are counts of labeled
records, not nonisomorphic graphs. The finite checks establish domain
coverage, matching, cover and size inequality. The subdivision argument
recovers every original length.

\Needspace{5\baselineskip}
\subsection{C.3. Two degree-three branches with degree-four
branches}\label{c.3.-two-degree-three-branches-with-degree-four-branches}

The number of degree-three vertices in each component is even. Under
\hyperlink{statement-lemma-4-1}{Lemma~4.1}, two such vertices must therefore lie in the same component.
The additional branch-degree profiles are \((3,3,4)\) and \((3,3,4,4)\).
Other components use some of the remaining degree-four vertices and are
covered by C.2.

Let \(l_i\) be loop counts and \(m_{ij}=m_{ji}\) nonloop multiplicities.
The full suppressed domain satisfies \[
2l_i+\sum_{j\ne i}m_{ij}=d_i,\quad
l_i,m_{ij}\in\mathbb Z_{\ge0},\quad
\text{connected support}.
\tag{C.1}\label{eq:C.1}
\] Each actual component gives a solution. Apply the parity/path rules
of C.2, keeping different paths internally disjoint. Loops of length one
or two are not introduced. Every original path is recovered by
two-vertex subdivisions. \hyperref[tab-mixed-branch]{Table 6} gives the
resulting domain counts.

\Needspace{10\baselineskip}
\begin{papertable}
\caption{Suppressed-graph and parity/path domains for the mixed branch-degree
profiles.}
\label{tab-mixed-branch}
\hypertarget{tab-mixed-branch}{}
\vspace{4pt}\centering\small\renewcommand{\arraystretch}{1.22}
\begin{tabular}{@{}>{\raggedright\arraybackslash}p{25.160mm}@{\hspace{3mm}}>{\raggedright\arraybackslash}p{25.160mm}@{\hspace{3mm}}>{\raggedright\arraybackslash}p{44.400mm}@{\hspace{3mm}}>{\raggedright\arraybackslash}p{25.160mm}@{\hspace{3mm}}>{\raggedright\arraybackslash}p{28.120mm}@{}}
\toprule
\textbf{Branch degrees} & \textbf{Labeled suppressed graphs} & \textbf{Topologies modulo degree-preserving branch permutations} & \textbf{Saved path records} & \textbf{Canonical parity/path records} \\
\midrule
\((3,3,4)\) & 7 & 5 & 339 & 249 \\
\((3,3,4,4)\) & 56 & 22 & 10,431 & 6,744 \\
\bottomrule
\end{tabular}
\end{papertable}

The reduced graphs have at most sixteen and twenty-two vertices,
respectively. Saved path counts are not counts of nonisomorphic simple
graphs. Each record has an explicit WKE witness.

One domain construction solves the degree equations recursively. A
separate checker enumerates off-diagonal multiplicity matrices, fixes
loop counts from the remaining even degrees, and tests connectivity. It
reconstructs maximal paths from each saved simple graph instead of
trusting its path description. Its domain modulo the stated symmetries
must equal the independently enumerated domain. Literal checks then
verify every witness. The predicates are thus both coverage and
validity, not just matching counts.

Lift the witnesses to the original paths by C.1 and combine components.
This proves \hyperlink{statement-lemma-4-1}{Lemma~4.1} subject to the finite certificate checks.
Suppression is not contraction of the ambient graph: the witness is
lifted back to its original link before \hyperlink{statement-lemma-2-3}{Lemma~2.3} is applied. \hfill\(\square\)

\Needspace{5\baselineskip}
\section{Appendix D. The quantitative surplus
count}\label{appendix-d.-the-quantitative-surplus-count}

Apply \hyperref[discharging-and-the-density-bound]{Section~3} in a nonempty reduction-free \(G\in\mathcal C\). Put
\(V_i=\{v:d_G(v)=i\}\), \(H=V_{\ge14}\), \(h=|H|\),
\(D=\sum_{v\in H}d_G(v)\), \(L=V_9\cup V_{10}\cup V_{11}\), and
\(\ell=|L|\). Write \(S_{12},S_{13},S_H\) for the surplus sums over the
indicated classes. Vertices in \(L\) have surplus zero, so
\(S=S_{12}+S_{13}+S_H\).

Each vertex in \(L\) has at least eight high neighbors and at most two
neighbors outside \(H\). Hence \[
8\ell\le D,\qquad e(V_{12},L)\le2\ell.
\tag{D.1}\label{eq:D.1}
\] For \(v\in V_{12}\), let \(r(v)=|N(v)\cap V_{\le12}|\le2\). Its
surplus is \((2-r(v))/20\). For \(v\in V_{13}\), the surplus is
\((10-|N(v)\cap V_{12}|)/20\). Thus \[
\sum_{v\in V_{12}}(2-r(v))=20S_{12},\qquad
|\{v\in V_{13}:|N(v)\cap V_{12}|<10\}|\le20S_{13}.
\tag{D.2}\label{eq:D.2}
\] For this count, call a vertex in \(V_{12}\) \emph{regular} if it has
exactly two degree-twelve neighbors and link \(T_A\). Let \(g\) count
these regular vertices, and \(B\) the other vertices in \(V_{12}\).
Vertices in \(V_{12}\) having fewer than two degree-twelve neighbors
number at most \[
\sum_{v\in V_{12}}(2-|N(v)\cap V_{12}|)
=20S_{12}+e(V_{12},L)\le20S_{12}+2\ell.
\tag{D.3}\label{eq:D.3}
\] Appendix A.4 leaves types \(T_C,T_D\) for other nonregular
degree-twelve vertices. A \(T_C\) vertex with two degree-twelve
neighbors forces both to have at most one neighbor of degree at most
twelve. Each such deficient neighbor is adjacent to at most one \(T_C\)
vertex, since its entire low neighborhood has size at most one.
Therefore \[
2\,|\{v\text{ of type }T_C\text{ with two twelve-neighbors}\}|
\le|\{v\in V_{12}:r(v)\le1\}|\le20S_{12}.
\tag{D.4}\label{eq:D.4}
\] Every \(T_D\) vertex has its label-11 neighbor in \(H\), so their
number is at most \(D\). Combining, \[
B\le30S_{12}+2\ell+D.
\tag{D.5}\label{eq:D.5}
\] A degree-thirteen vertex with ten degree-twelve neighbors must meet a
nonregular degree-twelve vertex. Otherwise all ten have type \(T_A\),
and the at-most-five versus at-least-seven edge count of \hyperref[appendix-b.-the-degree-twelvethirteen-endpoint]{Appendix~B}
applies without any global equality assumption. Every such nonregular
vertex has degree twelve, so, with \(n_{13}=|V_{13}|\), \[
n_{13}\le20S_{13}+12B.
\tag{D.6}\label{eq:D.6}
\] Every regular vertex has ten neighbors in \(V_{13}\cup H\), and every
degree-thirteen vertex has at most ten degree-twelve neighbors. Hence \[
10g\le10n_{13}+D.
\tag{D.7}\label{eq:D.7}
\] For \(N=|V(G)|=\ell+h+g+B+n_{13}\), equations \hyperref[eq:D.5]{(D.5)}--\hyperref[eq:D.7]{(D.7)} imply \[
\begin{aligned}
N&\le\ell+h+25B+40S_{13}+D/10\\
 &\le750S_{12}+40S_{13}+51\ell+h+(251/10)D\\
 &\le750S_{12}+40S_{13}+h+(1259/40)D.
\end{aligned}
\tag{D.8}\label{eq:D.8}
\] \Needspace{16\baselineskip}
The last step uses \hyperref[eq:D.1]{(D.1)}. \hyperref[discharging-and-the-density-bound]{Section~3} gives, for \(v\in H\), \[
s(v)\ge
\begin{cases}
1/30,&d_G(v)=14,\\
1/60,&d_G(v)=15,\\
2/15,&d_G(v)=16,\\
(79d_G(v)-1292)/100,&d_G(v)\ge17.
\end{cases}
\tag{D.9}\label{eq:D.9}
\] For each \(v\in H\), set \(d=d_G(v)\). These bounds give \[
1+(1259/40)d\le(56775/2)s(v).
\tag{D.10}\label{eq:D.10}
\] \Needspace{6\baselineskip}
At fifteen there is equality using the surplus floor; fourteen and
sixteen are direct rational checks. For \(d\ge17\), substitution gives
an affine difference with positive slope and positive value at
seventeen. Summing \hyperref[eq:D.10]{(D.10)}, then using \(750,40\le56775/2\) in \hyperref[eq:D.8]{(D.8)},
gives \hyperref[eq:5.2]{(5.2)}. Thus no equality hypothesis in \hyperlink{statement-theorem-1-1}{Theorem~1.1} is used in the
gap argument. \hfill\(\square\)

\Needspace{5\baselineskip}
\section{Appendix E. An ordinary-core sufficient
condition}\label{appendix-e.-an-ordinary-core-sufficient-condition}

Let \(Q_5,Q_9\) be the ordinary induced 5-core and 9-core of \(G\).
These are cores in \(G\), unlike the core in \(H_{\ge3}(G)\).

\par\addvspace{7pt}\Needspace{5\baselineskip}\noindent

\textbf{\statementtarget{statement-corollary-E-1}Corollary E.1.} If \(Q_5\in\mathcal D\), \(Q_9\in\mathcal C\),
and \(\operatorname{mad}(Q_9)\le25/2\), then \(\tau(G)\le2\nu(G)\). The
density condition may instead be \(\operatorname{mad}(Q_9)<25/2+\eta\).
No overlap restriction is imposed outside \(Q_5\).

\emph{Proof.} Every graph on at most four vertices is WKE. If its
maximum matching size is zero, it is edgeless. If it is one, its
nonempty component is a star or triangle, and one vertex covers all
edges outside a suitable matching edge. If it is two, a perfect matching
leaves a subgraph of a four-cycle, with a two-vertex cover.

Apply \hyperlink{statement-lemma-2-3}{Lemma~2.3} at current degree at most four, discarding isolates. The
survivor \(J_5\) is empty or has minimum degree at least five, so is
contained in the original \(Q_5\): every subgraph of minimum degree five
survives 5-core peeling. Its overlap conditions are inherited from
\(Q_5\). Appendix A.1 now reduces every center of current degree at most
eight. The survivor \(J_9\) is empty or has minimum degree nine, so lies
in the original \(Q_9\). Apply \hyperlink{statement-theorem-1-1}{Theorem~1.1} or 1.4 to
\(J_9\subseteq Q_9\), and telescope \hyperlink{statement-lemma-2-2}{Lemma~2.2} back to \(G\).

The survivors need not equal either full core. This uses hereditary
sufficient hypotheses on \(Q_9\), not an inequality comparing
\(\Phi(G)\) with \(\Phi(Q_9)\). \hfill\(\square\)

\Needspace{5\baselineskip}
\clearpage
\section{Use of AI tools}\label{use-of-ai-tools}

OpenAI tools, including Codex, were used to assist with the research and
preparation of this manuscript. The author takes full responsibility for
the content of this paper.

\section{Acknowledgements}\label{acknowledgements}

The Hifumi core is named after my son, Hifumi, whose name is written
{\hifumifont 一二三}, using the characters for one, two, and three. I also thank my
wife, Tomoko, for patiently tolerating my endless monologues about graph
theory at the dinner table.

\section{References}
\small\label{references}

\par\noindent\hangindent=1.8em\hangafter=1\hypertarget{bib-1}{}[1] Gregory J. Puleo. \emph{Tuza's Conjecture for Graphs with
Maximum Average Degree Less Than 7}. European Journal of Combinatorics
49 (2015), 134--152.
\href{https://doi.org/10.1016/j.ejc.2015.03.006}{doi:10.1016/j.ejc.2015.03.006}.
Passage numbering follows
\href{https://arxiv.org/abs/1308.2211v3}{arXiv:1308.2211v3}.

\par\noindent\hangindent=1.8em\hangafter=1\hypertarget{bib-2}{}[2] Amrit Kandasamy. \emph{Tuza's Conjecture for Graphs of Maximum
Degree at Most 10}.
\href{https://arxiv.org/abs/2610.15269v1}{arXiv:2610.15269v1}, 11
October 2026.

\par\noindent\hangindent=1.8em\hangafter=1\hypertarget{bib-3}{}[3] Amrit Kandasamy. \emph{Tuza at Maximum Average Degree at Most
45/4 and Book Number Three}. Supplementary document
outputs/tuza\_mad45\_4\_book3\_proof.txt in archive/verification.zip,
\href{https://github.com/Amkan2026/tuza-degree-ten}{tuza-degree-ten
repository}, snapshot retrieved 6 November 2026. Archive blob
\nolinkurl{7b0fddace492a587a9187e38d2dbb9317ba18bb4}; archive SHA-256
\nolinkurl{1c91ac287d56c46065c432fbc6a98835cfe50240b229ade07adc830f04429a5a}. This
is a supplementary statement, not Theorem 1.1 of [2].

\par\noindent\hangindent=1.8em\hangafter=1\hypertarget{bib-4}{}[4] Zhen Liu and Qinghou Zeng. \emph{Tuza's Conjecture under an
Edge-Count Condition}.
\href{https://arxiv.org/abs/2610.15539v1}{arXiv:2610.15539v1}, 11
October 2026.

\par\noindent\hangindent=1.8em\hangafter=1\hypertarget{bib-5}{}[5] Andrea Munaro. \emph{Triangle Packings and Transversals of Some
K4-Free Graphs}. Graphs and Combinatorics 34 (2018), 647--668.
\href{https://doi.org/10.1007/s00373-018-1903-y}{doi:10.1007/s00373-018-1903-y}.
\href{https://pureadmin.qub.ac.uk/ws/portalfiles/portal/186487135/manuscript.pdf}{Accepted
manuscript}.

\par\noindent\hangindent=1.8em\hangafter=1\hypertarget{bib-6}{}[6] Brendan D. McKay and Adolfo Piperno. \emph{Practical graph
isomorphism, II}. Journal of Symbolic Computation 60 (2014), 94--112.
\href{https://doi.org/10.1016/j.jsc.2013.09.003}{doi:10.1016/j.jsc.2013.09.003}.
\href{https://arxiv.org/abs/1301.1493}{arXiv:1301.1493}.

\end{document}